\documentclass{amsart}
\usepackage{amsfonts}

\usepackage{graphicx}
\usepackage{amscd}

\newtheorem{theorem}{Theorem}[section]
\theoremstyle{plain}
\newtheorem*{GMP}{Generalized Maximum Principle}
\newtheorem*{Theor.}{Theorem}

\newtheorem{lemma}[theorem]{Lemma}

\newtheorem{proposition}[theorem]{Proposition}

\theoremstyle{remark}
\newtheorem{remark}[theorem]{Remark}
\numberwithin{equation}{section}

\input{tcilatex}

\begin{document}
\title[Complete minimal hypersurfaces ]{Complete minimal hypersurfaces in the hyperbolic space $\mathbb{H}^{4}$ with
vanishing Gauss-Kronecker curvature}
\author{T. Hasanis}
\address[]{Department of Mathematics, University of Ioannina, 45110 Ioannina, Greece}
\email{thasanis@cc.uoi.gr, me00499@cc.uoi.gr, tvlachos@cc.uoi.gr}
\author{A. Savas-Halilaj}
\author{T. Vlachos}

\begin{abstract}
We investigate 3-dimensional complete minimal hypersurfaces in the
hyperbolic space $\mathbb{H}^{4}$ with Gauss-Kronecker curvature identically
zero. More precisely, we give a classification of complete minimal
hypersurfaces with Gauss-Kronecker curvature identically zero, nowhere
vanishing second fundamental form and scalar curvature bounded from below.
\end{abstract}

\maketitle

\section{Introduction}

In order to study the rigidity of minimal hypersurfaces, Dajczer and Gromoll 
\cite{DG} invented the so called \textit{Gauss parametrization}. As a
by-product of this approach they were able to describe locally the minimal
hypersurfaces in the $\left( n+1\right) $-dimensional space form whenever
the rank of the nullity distribution is constant. Almeida and Brito \cite{AB}
initiated the study of compact minimal hypersurfaces in the unit sphere $%
\mathbb{S}^{4}$ with vanishing Gauss-Kronecker curvature. In fact they
proved that such compact hypersurfaces are boundaries of tubes of minimal
2-spheres in $\mathbb{S}^{4}$, provided that the second fundamental form
never vanishes. Ramanathan \cite{R} extended this result and allowed points
where the second fundamental form is zero. In \cite{H1}, \cite{H2}\ the
authors extended the above results to complete minimal hypersurfaces in the
Euclidean space $\mathbb{R}^{4}$ or in the unit sphere $\mathbb{S}^{4}$.

The aim of this paper is to study complete minimal hypersurfaces in the
4-dimensional hyperbolic space $\mathbb{H}^{4}$ with identically zero
Gauss-Kronecker curvature. We remind that the Gauss-Kronecker curvature is
the product of the principal curvatures. In fact, we deal with minimal
hypersurfaces whose second fundamental form is nowhere zero, which is
equivalent to the fact that the nullity distribution is one dimensional.

It turns out that such hypersurfaces are closely related to stationary
spacelike surfaces in the de Sitter space $\mathbb{S}_{1}^{4}$, which is the
Lorentzian unit sphere in the flat Lorentzian space $\mathbb{R}_{1}^{5}$.
More precisely, Dajczer and Gromoll \cite{DG} noticed that the unit normal
bundle of a stationary spacelike surface in $\mathbb{S}_{1}^{4}$ gives rise
to a minimal hypersurface in $\mathbb{H}^{4}$ with vanishing Gauss-Kronecker
curvature and nowhere zero second fundamental form (for details see Section
2) and, conversely, any such hypersurface is obtained, at least locally, in
this way.

We are interested in the classification of complete minimal hypersurfaces in 
$\mathbb{H}^{4}$ with Gauss-Kronecker curvature identically zero. At first
we show that there exists an abundance of such hypersurfaces. To this
purpose we focus on a class of stationary spacelike surfaces in $\mathbb{S}%
_{1}^{4}$, namely those with vanishing normal curvature. We establish a
correspondence between stationary spacelike surfaces in $\mathbb{S}_{1}^{4}$
with identically zero normal curvature and minimal surfaces in umbilical
hypersurfaces in\ the hyperbolic space $\mathbb{H}^{4}$.

Our examples of complete minimal hypersurfaces in $\mathbb{H}^{4}$ with
vanishing Gauss-Kronecker curvature arise as suspensions of complete minimal
surfaces in horospheres or equidistant hypersurfaces of $\mathbb{H}^{4}$.
Then the question whether these are the only examples comes naturally.

We prove that these suspensions are in fact the only complete minimal
hypersurfaces in $\mathbb{H}^{4}$ with identically zero Gauss-Kronecker
curvature under the assumptions that the second fundamental form is nowhere
zero and the scalar curvature is bounded from below.

The paper is organized as follows: In Section 2 we study stationary
spacelike surfaces in the de Sitter space $\mathbb{S}_{1}^{4}$, we define
the polar map and show that it induces minimal hypersurfaces in $\mathbb{H}%
^{4}$ with vanishing Gauss-Kronecker curvature. In particular, we prove some
auxiliary results about stationary spacelike surfaces in $\mathbb{S}_{1}^{4}$
with identically zero normal curvature. Moreover, we furnish a method to
produce complete minimal hypersurfaces in $\mathbb{H}^{4}$ with identically
zero Gauss-Kronecker curvature. Section 3 is devoted to the local theory of
minimal hypersurfaces in $\mathbb{H}^{4}$ with vanishing Gauss-Kronecker
curvature. Furthermore, we give some auxiliary results. Finally, in Section
4 we state and prove the main result of this paper.

\section{The polar map of stationary surfaces in the de Sitter space}

At first we set up our notation. Denote by $\mathbb{R}_{1}^{5}$ the real
vector space $\mathbb{R}^{5}$ endowed with the Lorentzian metric tensor $%
\left\langle \ ,\ \right\rangle $ given by 
\begin{equation*}
\left\langle x,y\right\rangle =-x_{0}y_{0}+\dsum_{i=1}^{4}x_{i}y_{i},
\end{equation*}
where $x=\left( x_{0},x_{1},x_{2},x_{3},x_{4}\right) $, $y=\left(
y_{0},y_{1},y_{2},y_{3},y_{4}\right) \in \mathbb{R}^{5}$. We shall use the
Minkowski model for the simply connected \textit{hyperbolic space} of
constant sectional curvature $-1$, which is the hyperquadric 
\begin{equation*}
\mathbb{H}^{4}=\left\{ x\in \mathbb{R}_{1}^{5}:\left\langle x,x\right\rangle
=-1,x_{0}>0\right\} ,
\end{equation*}
Moreover, the hyperquadric 
\begin{equation*}
\mathbb{S}_{1}^{4}=\left\{ x\in \mathbb{R}_{1}^{5}:\left\langle
x,x\right\rangle =1\right\} ,
\end{equation*}
is the standard model for the simply connected Lorentzian space form of
constant curvature 1, and is called the\ \textit{de Sitter space}.

Consider a 2-dimensional manifold $M^{2}$. An immersion $g:M^{2}\rightarrow 
\mathbb{S}_{1}^{4}$ is called \textit{spacelike }if the induced metric on $%
M^{2}$ via $g$ is Riemannian, which as usual will be denoted again by $%
\left\langle \ ,\ \right\rangle $. Let $i:\mathbb{S}_{1}^{4}\rightarrow 
\mathbb{R}_{1}^{5}$ be the inclusion map. Denote by 
\begin{eqnarray*}
\left( i\circ g\right) ^{\ast }\left( T\mathbb{R}_{1}^{5}\right) &=&\left\{
\left( x,w\right) :x\in M^{2},w\in T_{g\left( x\right) }\mathbb{R}%
_{1}^{5}\right\} , \\
g^{\ast }\left( T\mathbb{S}_{1}^{4}\right) &=&\left\{ \left( x,w\right)
:x\in M^{2},w\in T_{g\left( x\right) }\mathbb{S}_{1}^{4}\right\} ,
\end{eqnarray*}
the \textit{induced bundles} of $i\circ g$ and $g$, respectively. The 
\textit{normal bundle} $\mathcal{N}\left( g\right) $ of $g$ is given by 
\begin{equation*}
\mathcal{N}\left( g\right) =\left\{ \left( x,w\right) \in g^{\ast }\left( T%
\mathbb{S}_{1}^{4}\right) :w\perp dg\left( T_{x}M^{2}\right) \right\} .
\end{equation*}
Also we denote by $\overline{\nabla }$, $\overset{g}{\nabla }$ the \textit{%
connections of the induced bundles} of $i\circ g$ and $g$, respectively, and
by $\overset{g}{D}$ the \textit{connection of the normal bundle} of $g$.
Given a normal vector field $\eta $ along $g$ and a tangent vector $X$ of $%
M^{2}$, then we have 
\begin{equation*}
\overline{\nabla }_{X}\eta =\overset{g}{\nabla }_{X}\eta -\left\langle
dg\left( X\right) ,\eta \right\rangle g.
\end{equation*}
The \textit{second fundamental} \textit{form} $II$ of $g$, is given by the 
\textit{Gauss formula} 
\begin{equation*}
II\left( X,Y\right) =\overset{g}{\nabla }_{X}dg\left( Y\right) -dg\left(
\nabla _{X}Y\right) ,
\end{equation*}
where $X$,$Y$ are tangent vector fields of $M^{2}$ and $\nabla $ stands for
the Levi-Civita connection of the induced metric\ on $M^{2}$. The
selfadjoint operator $A_{\eta }$ defined by 
\begin{equation*}
\left\langle A_{\eta }X,Y\right\rangle =\left\langle II\left( X,Y\right)
,\eta \right\rangle
\end{equation*}
is called the \textit{shape operator} of $g$ relative to $\eta $. Moreover,
the \textit{Weingarten formula} is 
\begin{equation*}
\overset{g}{\nabla }_{X}\eta =-dg\left( A_{\eta }X\right) +\overset{g}{D}%
_{X}\eta \text{.}
\end{equation*}
A point $x\in M^{2}$ is called \textit{totally geodesic point of }$g$ if and
only if $II_{x}=0$. If each point of $M^{2}$ is a totally geodesic point of $%
g$, then $g$ is called \textit{totally geodesic immersion}. Since $\mathbb{S}%
_{1}^{4}$ is a Lorentz manifold of constant sectional curvature, one can
derive the \textit{equations of Codazzi and Ricci}, which are respectively\ 
\begin{eqnarray*}
\left( \nabla _{X}A_{\eta }\right) Y+A_{\overset{g}{D}_{Y}\eta }X &=&\left(
\nabla _{Y}A_{\eta }\right) X+A_{\overset{g}{D}_{X}\eta }Y, \\
\left\langle R^{D}\left( X,Y\right) \eta _{1},\eta _{2}\right\rangle
&=&\left\langle [A_{\eta _{1}},A_{\eta _{2}}]X,Y\right\rangle ,
\end{eqnarray*}
where $X$,$Y$ are tangent vector fields of $M^{2}$, $R^{D}$ is the curvature
tensor of $\overset{g}{D}$ and $\eta $, $\eta _{1}$, $\eta _{2}$ are normal
vector fields along $g$.

Consider now an orthonormal adapted frame field $\left\{
e_{1},e_{2};e_{3},e_{4}\right\} $ along $g$, where $e_{4}$ is timelike.
Then, we have 
\begin{equation*}
II\left( X,Y\right) =\left\langle A_{3}X,Y\right\rangle e_{3}-\left\langle
A_{4}X,Y\right\rangle e_{4},
\end{equation*}
where $A_{3}$, $A_{4}$ are the shape operators of $g$ with respect to the
directions $e_{3}$ and $e_{4}$. The \textit{mean curvature vector field} $H$
is given by 
\begin{equation*}
H=\frac{1}{2}\left( traceA_{3}\right) e_{3}-\frac{1}{2}\left(
traceA_{4}\right) e_{4}.
\end{equation*}
The immersion $g$ is called \textit{stationary,} whenever $H\equiv 0$. The 
\textit{Gaussian curvature} $K$ of the induced metric\ is 
\begin{equation*}
K=1+\det A_{3}-\det A_{4}.
\end{equation*}
We denote by $\left\{ \omega _{1},\omega _{2}\right\} $ the \textit{dual
frame }of $\left\{ e_{1},e_{2}\right\} $ and by $\omega _{34}$ the \textit{%
connection form of the normal bundle of} $g$, which is determined by 
\begin{equation*}
\omega _{34}\left( X\right) =-\left\langle \overset{g}{D}_{X}e_{3},e_{4}%
\right\rangle .
\end{equation*}
The \textit{normal curvature} $K^{\bot }$ of $g$ is given by 
\begin{equation*}
K^{\bot }=\left\langle R^{D}\left( e_{1},e_{2}\right)
e_{3},e_{4}\right\rangle =\left\langle \left[ A_{3,}A_{4}\right]
e_{1},e_{2}\right\rangle .
\end{equation*}
We recall that 
\begin{equation}
d\omega _{34}=-K^{\perp }\omega _{1}\wedge \omega _{2}.
\end{equation}

Assume now that $M^{2}$ is an oriented, 2-dimensional Riemannian manifold
and $g:M^{2}\rightarrow \mathbb{S}_{1}^{4}$ is a stationary isometric
immersion. It is well known (cf. \cite{AR}) that there exists a holomorphic
quadric differential on $M^{2}$, the so called \textit{Hopf differential}. A
point $x\in M^{2}$ is a zero of the Hopf differential if and only if $%
K\left( x\right) =1$ and $K^{\bot }\left( x\right) =0$. Such a point is
called a \textit{superminimal point of }$g$. The immersion $g$ is called 
\textit{superminimal} if each point of $M^{2}$ is a superminimal point of $g$%
. From the holomorphicity of the Hopf differential, it follows that either $%
g $ is superminimal or the superminimal points are isolated.

Consider the \textit{unit normal bundle} $\mathcal{N}^{1}\left( g\right) $\
of $g$, defined by 
\begin{equation*}
\mathcal{N}^{1}\left( g\right) =\left\{ \left( x,w\right) \in \mathcal{N}%
\left( g\right) :\left\langle w,w\right\rangle =-1\right\} .
\end{equation*}
Denote by $\pi :\mathcal{N}^{1}\left( g\right) \rightarrow M^{2}$ the
projection to the first factor\ and by\ $\Psi _{g}:\mathcal{N}^{1}\left(
g\right) \rightarrow \mathbb{H}^{4}$ the projection to the second factor.
The map $\Psi _{g}$ is called the \textit{polar map} associated with $g$.

Throughout the paper we follow the above mentioned notation and assume that
all manifolds under consideration are connected, unless otherwise stated.

The following proposition was essentially proved by Dajczer and Gromoll in 
\cite{DG}. In order to make the paper self-contained we shall give here a
proof that fits our exposition. This proposition furnishes a method for
producing minimal hypersurfaces in $\mathbb{H}^{4}$ with Gauss-Kronecker
curvature identically zero and establishes the close relation between them
and stationary surfaces in $\mathbb{S}_{1}^{4}$.

\begin{proposition}
\label{p1}Let $M^{2}$ be a $2$-dimensional Riemannian manifold and $%
g:M^{2}\rightarrow \mathbb{S}_{1}^{4}$ a stationary isometric immersion. Then

\begin{enumerate}
\item[$\left( i\right) $]  The polar map $\Psi _{g}$ associated with $g$ is
regular at $\left( y,w\right) \in \mathcal{N}^{1}\left( g\right) $ if and
only if the second fundamental form of $g$ is non-singular in the direction $%
w$.

\item[$\left( ii\right) $]  On the open set of its regular points, $\Psi
_{g} $ is a minimal immersion in $\mathbb{H}^{4}$ with Gauss-Kronecker
curvature identically zero and nowhere vanishing second fundamental form.

\item[$\left( iii\right) $]  If $x$ is a point on $M^{2}$ where the normal
curvature of $g$ is not zero, then $\Psi _{g}$ is regular on the fiber of $%
\mathcal{N}^{1}\left( g\right) $ over $x$. Furthermore, if $x$ is not a
totally geodesic point of $g$ and $K\left( x\right) \geq 1$, then $\Psi _{g}$
is regular on the fiber of $\mathcal{N}^{1}\left( g\right) $ over $x$.
\end{enumerate}
\end{proposition}

\begin{proof}
Choose an adapted orthonormal frame field $\left\{
e_{1},e_{2};e_{3},e_{4}\right\} $, along $g$ defined on an open set $%
U\subset M^{2}$, where $e_{4}$ is timelike. We parametrize $\pi ^{-1}\left(
U\right) $ by $U\times \mathbb{R}$ via the map 
\begin{equation*}
\left( x,t\right) \rightarrow \left( x,\sinh te_{3}\left( x\right) +\cosh
te_{4}\left( x\right) \right) .
\end{equation*}
Then $\Psi _{g}\left( x,t\right) =\sinh te_{3}\left( x\right) +\cosh
te_{4}\left( x\right) $. For the sake of convenience we set $W\left(
x,t\right) =\sinh te_{3}\left( x\right) +\cosh te_{4}\left( x\right) $.
Obviously $w=W\left( y,t_{0}\right) $, for some $t_{0}\in \mathbb{R}$.

$\left( i\right) $ Calculating the differential of $\Psi _{g}$ at the point $%
\left( y,w\right) $ we have, 
\begin{equation*}
d\Psi _{g}\left( \frac{\partial }{\partial t}\right) =\cosh t_{0}e_{3}\left(
y\right) +\sinh t_{0}e_{4}\left( y\right) ,
\end{equation*}
and 
\begin{eqnarray}
d\Psi _{g}\left( X\right) &=&\overline{\nabla }_{X}W=\overset{g}{\nabla }%
_{X}W \\
&=&-dg\left( A_{w}X\right) +\overset{g}{D}_{X}W  \notag \\
&=&-dg\left( A_{w}X\right) +\omega _{34}\left( X\right) \left( \cosh
t_{0}e_{3}\left( y\right) +\sinh t_{0}e_{4}\left( y\right) \right)  \notag \\
&=&-dg\left( A_{w}X\right) +\omega _{34}\left( X\right) d\Psi _{g}\left( 
\frac{\partial }{\partial t}\right) ,  \notag
\end{eqnarray}
for each $X\in T_{y}M^{2}$. From the above relations it follows that $\left(
y,w\right) $ is a regular point of $\Psi _{g}$ if and only if $\det
A_{w}\left( y\right) \neq 0$.

$\left( ii\right) $ The vector field $\xi $ given by $\xi \left( x,t\right)
=g\left( x\right) $, $\left( x,t\right) \in U\times \mathbb{R}$, is a\ unit
normal vector field along $\Psi _{g}$. Denote by $A_{\xi }$ the
corresponding shape operator and by $\pi _{1}:U\times \mathbb{R\rightarrow }%
U $, $\pi _{2}:U\times \mathbb{R\rightarrow R}$ the corresponding projection
maps.\ Using the Weingarten formula, we have 
\begin{equation*}
0=d\xi \left( \frac{\partial }{\partial t}\right) =-d\Psi _{g}\left( A_{\xi }%
\frac{\partial }{\partial t}\right) .
\end{equation*}
Consequently, the\ Gauss-Kronecker curvature of $\Psi _{g}$ is identically
zero. Moreover using $\left( 2.2\right) $, we get 
\begin{eqnarray*}
-dg\left( X\right) &=&-d\xi \left( X\right) =d\Psi _{g}\left( A_{\xi
}X\right) \\
&=&d\Psi _{g}\left( d\pi _{1}\left( A_{\xi }X\right) +d\pi _{2}\left( A_{\xi
}X\right) \right) \\
&=&-dg\left( A_{w}\left( d\pi _{1}\left( A_{\xi }X\right) \right) \right) +%
\overset{g}{D}_{d\pi _{1}\left( A_{\xi }X\right) }W+d\Psi _{g}\left( d\pi
_{2}\left( A_{\xi }X\right) \right) ,
\end{eqnarray*}
for each $X\in T_{y}M^{2}$. So 
\begin{equation*}
d\pi _{1}\left( A_{\xi }X\right) =A_{w}^{-1}X,
\end{equation*}
and $\Psi _{g}$ has principal curvatures 
\begin{equation*}
k_{1}\left( x,w\right) =-k_{3}\left( x,w\right) =\frac{1}{\sqrt{-\det
A_{w}\left( x\right) }},\ k_{2}\left( x,w\right) =0\text{.}
\end{equation*}

$\left( iii\right) $ Suppose that $K^{\bot }\left( x\right) \neq 0$. Then,
obviously$\ \det A_{w}\left( x\right) \neq 0$, for each $w$ on the fiber of $%
\mathcal{N}^{1}\left( g\right) $ over $x$. Assume now that $x$ is not a
totally geodesic point of $g$, $K\left( x\right) \geq 1$ and that there
exists a vector$\ w$ on the fiber of $\mathcal{N}^{1}\left( g\right) $ over $%
x$ such that$\ \det A_{w}\left( x\right) =0$. Let $\eta $ be a unit normal
vector in the normal bundle of $g$ such that $\left\langle \eta
,w\right\rangle =0$. Then, because $x$ is not a totally geodesic point, it
follows that $K\left( x\right) =1+\det A_{\eta }\left( x\right) -\det
A_{w}\left( x\right) =1+\det A_{\eta }\left( x\right) <1$, which is a
contradiction. This completes the proof.
\end{proof}

\begin{remark}
We shall see in Section 3 that every minimal hypersurface in the hyperbolic
space $\mathbb{H}^{4}$ with Gauss-Kronecker curvature identically zero and
nowhere vanishing second fundamental form can be obtained, at least locally,
as in Proposition \ref{p1}$\left( ii\right) $.
\end{remark}

We focus now on the class of stationary spacelike minimal surfaces in $%
\mathbb{S}_{1}^{4}$ with identically zero normal curvature. This class will
play a crucial role in the classification of complete minimal hypersurfaces
in the hyperbolic space with zero Gauss-Kronecker curvature. We remind here
that the\ totally geodesic submanifolds of $\mathbb{S}_{1}^{4}$ arise as
intersections of $\mathbb{S}_{1}^{4}$ with linear subspaces of $\mathbb{R}%
_{1}^{5}$. The following result is due to Alias and Palmer \cite{AR}. For
the sake of completeness we give another short proof.

\begin{proposition}
\label{p2}Let $M^{2}$ be a 2-dimensional Riemannian manifold and $%
g:M^{2}\rightarrow \mathbb{S}_{1}^{4}$ be a stationary isometric immersion.
Then $K^{\bot }\equiv 0$ if and only if $g\left( M^{2}\right) $ is contained
in a totally geodesic hypersurface $L^{3}$ of$\ \mathbb{S}_{1}^{4}$, i.e.,
there exists a vector $w$ such that $\left\langle g\left( x\right)
,w\right\rangle =0$, for each $x\in M^{2}$. Moreover,

\begin{enumerate}
\item[$\left( i\right) $]  $w$ is spacelike if and only if $K\geq 1$ and $%
K\not\equiv 1$,

\item[$\left( ii\right) $]  $w$ is timelike if and only if $K\leq 1$ and $%
K\not\equiv 1$,

\item[$\left( iii\right) $]  $w$ may be chosen to be null if and only if $%
K\equiv 1$.
\end{enumerate}
\end{proposition}

\begin{proof}
Denote by $M_{1}$ the set of superminimal points of $g$. We distinguish two
cases.

\textbf{Case 1. }Assume that $M_{1}$ consists of isolated points only.%
\textbf{\ }Consider a non superminimal point $x$ of $g$ and let $\left\{
e_{1},e_{2};e_{3},e_{4}\right\} $ be an adapted\ orthonormal frame field
defined on an open set $U_{x}\subset M^{2}-M_{1}$ around $x$, $e_{4}$ being
timelike. We may suppose that the shape operators associated with $e_{3}$
and $e_{4}$ are 
\begin{equation*}
A_{3}\sim \left( 
\begin{array}{cc}
\kappa & \text{ \ }0 \\ 
0 & -\kappa
\end{array}
\right) ,\text{ \ }A_{4}\sim \left( 
\begin{array}{cc}
\mu & \text{ \ }0 \\ 
0 & -\mu
\end{array}
\right) .
\end{equation*}
Because $K\neq 1$ on $U_{x},$ it follows that $\kappa ^{2}\neq \mu ^{2}$.
Define now the vector fields 
\begin{equation*}
\overline{e}_{3}=\frac{1}{\sqrt{\left| \mu ^{2}-\kappa ^{2}\right| }}\left(
\mu e_{3}-\kappa e_{4}\right) ,\text{ }\overline{e}_{4}=\frac{1}{\sqrt{%
\left| \mu ^{2}-\kappa ^{2}\right| }}\left( \kappa e_{3}-\mu e_{4}\right) .
\end{equation*}
Then the shape operators corresponding to the directions $\overline{e}_{3}$
and $\overline{e}_{4}$ are 
\begin{equation*}
\overline{A}_{3}=0\text{ \ and \ }\overline{A}_{4}\sim \left( 
\begin{array}{cc}
\frac{\kappa ^{2}-\mu ^{2}}{\sqrt{\left| \mu ^{2}-\kappa ^{2}\right| }} & 0
\\ 
0 & -\frac{\kappa ^{2}-\mu ^{2}}{\sqrt{\left| \mu ^{2}-\kappa ^{2}\right| }}
\end{array}
\right) .
\end{equation*}
From the Codazzi equation we get, 
\begin{equation*}
A_{\overset{g}{D}_{e_{1}}\overline{e}_{3}}e_{2}=A_{\overset{g}{D}_{e_{2}}%
\overline{e}_{3}}e_{1},
\end{equation*}
or, equivalently, 
\begin{equation*}
\overline{\omega }_{34}\left( e_{1}\right) \left( \overline{A}%
_{4}e_{2}\right) =\overline{\omega }_{34}\left( e_{2}\right) \left( 
\overline{A}_{4}e_{1}\right) ,
\end{equation*}
where $\overline{\omega }_{34}$ stands for the connection form on the normal
bundle of $g$ with respect to the frame $\left\{ \overline{e}_{3},\overline{e%
}_{4}\right\} $. Thus $\overline{\omega }_{34}=0$, and so the vector field $%
w=\overline{e}_{3}$ is constant along $g$ and $\left\langle g,w\right\rangle
=0$ on $U_{x}$. This means that $g\left( U_{x}\right) $ is contained in a
totally geodesic hypersurface $L^{3}$ of $\mathbb{S}_{1}^{4}$. Note that the
normal vector field $w$ satisfies 
\begin{equation*}
\left\langle w,w\right\rangle =\frac{\mu ^{2}-\kappa ^{2}}{\left| \mu
^{2}-\kappa ^{2}\right| }=\frac{K-1}{\left| K-1\right| }.
\end{equation*}
Suppose now that $U_{y}\subset M^{2}-M_{1}$ is an open set around another
point $y\in M^{2}-M_{1}$, such that $U_{x}\cap U_{y}\neq \emptyset $ and $%
g\left( U_{y}\right) \subset \overline{L}^{3}$, where $\overline{L}^{3}$ is
a totally geodesic hypersurface of $\mathbb{S}_{1}^{4}$, with normal vector $%
\overline{w}$. We claim that $\overline{L}^{3}=L^{3}$. To this purpose, it
is enough to prove that $w$ and $\overline{w}$ are linearly dependent.
Suppose in the contrary that these are linearly independent. Then $g$ is
totally geodesic on $U_{x}\cap U_{y}$, a contradiction since $K\neq 1$ on $%
U_{x}\cap U_{y}$.\ So we deduce that $g\left( M^{2}\right) $ is contained in
a totally geodesic hypersurface $L^{3}$ of $\mathbb{S}_{1}^{4}$, whose
normal vector satisfies 
\begin{equation*}
\left\langle w,w\right\rangle =\frac{K-1}{\left| K-1\right| },
\end{equation*}
on $M^{2}-M_{1}$. This completes the proof of parts $\left( i\right) $ and $%
\left( ii\right) $.

\textbf{Case 2.} Suppose that $M_{1}=M^{2}$. Then $g$ is superminimal and $%
K\equiv 1$. Because $K^{\bot }\equiv 0$, around each point we may choose a
parallel orthonormal frame field $\left\{ \eta _{1},\eta _{2}\right\} $ in
the normal bundle of $g$. Let $\left( u,v\right) $ be local isothermal
coordinates. The complex valued functions, 
\begin{equation*}
\sigma _{i}\left( u,v\right) =\left\langle A_{\eta _{i}}\frac{\partial }{%
\partial u},\frac{\partial }{\partial u}\right\rangle -\sqrt{-1}\left\langle
A_{\eta _{i}}\frac{\partial }{\partial u},\frac{\partial }{\partial v}%
\right\rangle ,\text{ \ }i=1,2,
\end{equation*}
are holomorphic and their zeroes are precisely the totally geodesic points
of $g$. Hence the set $M_{0}$ of totally geodesic points of $g$ either
coincides with $M^{2}$ or consists of isolated points. In the case where $%
M_{0}=M^{2}$, $g\left( M^{2}\right) $ is contained in a appropriate totally
geodesic hypersurface of $\mathbb{S}_{1}^{4}$ whose normal vector is null.
Suppose now that $M_{0}$ consists of isolated points. Then the set $%
M^{2}-M_{0}$ is open and connected. Consider a non totally geodesic point $%
x\in M^{2}$ and let $\left\{ e_{1},e_{2};e_{3},e_{4}\right\} $ be an
adapted\ orthonormal frame field defined on a simply connected neighborhood $%
U\subset M^{2}-M_{0}$ of $x$, $e_{4}$ being timelike. We may suppose that $%
A_{3}=A_{4}$. Since $K^{\bot }\equiv 0$, from $\left( 2.1\right) $ it
follows that $d\omega _{34}=0$. Thus there exists a smooth function $\theta $
on $U$ such that $\omega _{34}=d\theta $. Define the vector field 
\begin{equation*}
w=e^{\theta }\left( e_{3}-e_{4}\right) .
\end{equation*}
For each tangent vector field $X$ of $U$, we have 
\begin{eqnarray*}
dw\left( X\right) &=&e^{\theta }X\left( \theta \right) \left(
e_{3}-e_{4}\right) +e^{\theta }\overset{g}{\nabla }_{X}\left(
e_{3}-e_{4}\right) \\
&=&e^{\theta }X\left( \theta \right) \left( e_{3}-e_{4}\right) +e^{\theta
}\omega _{34}\left( X\right) \left( e_{4}-e_{3}\right) \\
&=&0.
\end{eqnarray*}
Therefore $w$ is constant and $\left\langle g,w\right\rangle =0$. Thus, $%
g\left( U\right) $ is contained in a totally geodesic hypersurface $L^{3}$
of $\mathbb{S}_{1}^{4}$ whose normal vector $w$ satisfies $\left\langle
w,w\right\rangle =0$. Arguing as in Case 1 we can prove that $g\left(
M^{2}\right) \subset L^{3}$.\ This completes the proof.
\end{proof}

The following proposition provides a way to produce all spacelike stationary
surfaces in $\mathbb{S}_{1}^{4}$ with normal curvature identically zero. For
the sake of convenience, we introduce the following notation. Let $%
h:M^{2}\rightarrow Q^{3}$ be an isometric immersion, where $M^{2}$ is a $2$%
-dimensional, oriented Riemannian manifold and $Q^{3}$ an umbilical
hypersurface of $\mathbb{H}^{4}$. The orientation $N$ of $h$ gives rise in a
natural way to a map $\widehat{h}:=N:M^{2}\rightarrow \mathbb{S}_{1}^{4}$
which is called\textit{\ }the \textit{associate} of $h$.

We recall here that the umbilical hypersurfaces of $\mathbb{H}^{4}$ arise as
intersections of $\mathbb{H}^{4}$ with affine hyperplanes of $\mathbb{R}%
_{1}^{5}$. Moreover, an umbilical hypersurface $Q^{3}$ of $\mathbb{H}^{4}$
has positive, negative or zero sectional curvature if $Q^{3}$ is a \textit{%
geodesic sphere}, an \textit{equidistant hypersurface} or a \textit{%
horosphere}, respectively.

\begin{proposition}
\label{p3}Let $M^{2}$ be a 2-dimensional, oriented Riemannian manifold.

\begin{enumerate}
\item[$\left( i\right) $]  If $h:M^{2}\rightarrow Q^{3}$ is\ a minimal
isometric immersion without totally geodesic points, where $Q^{3}$ is an
umbilical hypersurface of $\mathbb{H}^{4}$, then its associate $\widehat{h}%
:M^{2}\rightarrow \mathbb{S}_{1}^{4}$ is a spacelike stationary immersion
with normal curvature identically zero without totally geodesic points.

\item[$\left( ii\right) $]  Conversely, assume that $g:M^{2}\rightarrow 
\mathbb{S}_{1}^{4}$ is a stationary isometric immersion with normal
curvature identically zero without totally geodesic points. Then there exist
a vector $w$ and a totally geodesic point free minimal immersion $%
h:M^{2}\rightarrow Q^{3}$, where $Q^{3}$ is an umbilical hypersurface of $%
\mathbb{H}^{4}$, with sectional curvature $K_{Q^{3}}=-\left\langle
w,w\right\rangle $ such that $\left\langle g\left( x\right) ,w\right\rangle
=0$, for each $x\in M^{2}$ and $g=\widehat{h}$.
\end{enumerate}
\end{proposition}

\begin{proof}
$\left( i\right) $ Let $\eta $ be a unit normal vector of $Q^{3}$ in $%
\mathbb{H}^{4}$ and $A_{\eta }$ the corresponding shape operator. Obviously, 
$A_{\eta }=\alpha I$, for some $\alpha \in \mathbb{R}$. We denote by $N$ the
orientation of $h$ and by $A_{N}$ the corresponding shape operator. For each
tangent vector $X$ of $M^{2}$, we have 
\begin{equation}
d\widehat{h}\left( X\right) =-dh\left( A_{N}X\right) .
\end{equation}
Hence $\widehat{h}$ is an immersion and the metric $\left\langle
X,Y\right\rangle _{\widehat{h}}=\left\langle A_{N}^{2}X,Y\right\rangle $
induced by $\widehat{h}$ on $M^{2}$ is Riemannian, where $\left\langle \text{
},\text{ }\right\rangle $ stands for the Riemannian metric of $M^{2}$.
Moreover, the vector fields $\left\{ e_{3}:=\eta \circ h,e_{4}:=h\right\} $
constitute a frame field in the normal bundle of $\widehat{h}$. Denote by $%
\widehat{A}_{3}$, $\widehat{A}_{4}$ the shape operators of $\widehat{h}$
with respect to directions $e_{3}$ and $e_{4}$. Using the Weingarten formula
and $\left( 2.3\right) $, we get 
\begin{equation*}
de_{3}\left( X\right) =\overset{\widehat{h}}{\nabla }_{X}e_{3}=-d\widehat{h}%
\left( \widehat{A}_{3}X\right) +\overset{\widehat{h}}{D}_{X}e_{3}=dh\left(
A_{N}\widehat{A}_{3}X\right) +\overset{\widehat{h}}{D}_{X}e_{3}.
\end{equation*}
Since, 
\begin{equation*}
de_{3}\left( X\right) =d\eta \left( dh\left( X\right) \right) =-A_{\eta
}\left( dh\left( X\right) \right) =-\alpha dh\left( X\right) ,
\end{equation*}
it follows that 
\begin{equation}
\widehat{A}_{3}=-\alpha A_{N}^{-1}\text{ \ and \ }\overset{\widehat{h}}{D}%
_{X}e_{3}=0.
\end{equation}
Moreover, 
\begin{equation*}
dh\left( X\right) =\overset{\widehat{h}}{\nabla }_{X}e_{4}=-d\widehat{h}%
\left( \widehat{A}_{4}X\right) +\overset{\widehat{h}}{D}_{X}e_{4}=dh\left(
A_{N}\widehat{A}_{4}X\right) +\overset{\widehat{h}}{D}_{X}e_{4}.
\end{equation*}
Therefore 
\begin{equation}
\widehat{A}_{4}=A_{N}^{-1}\text{ \ and \ }\overset{\widehat{h}}{D}%
_{X}e_{4}=0.
\end{equation}
From $\left( 2.4\right) $ and $\left( 2.5\right) $ we deduce that $\widehat{h%
}$ is stationary with normal curvature identically zero.

$\left( ii\right) $ Suppose now that $g:M^{2}\rightarrow \mathbb{S}_{1}^{4}$
is a stationary isometric immersion with normal curvature identically zero
without totally geodesic points. According to Proposition \ref{p2}, $g\left(
M^{2}\right) $ is contained in a totally geodesic hypersurface of $\mathbb{S}%
_{1}^{4}$ with normal vector $w$. Without loss of generality, we may assume
that $\left\langle w,w\right\rangle \leq 1$. We distinguish two cases.

\textbf{Case 1. }Assume that $w$ is not null. We set $a:=\sqrt{%
1-\left\langle w,w\right\rangle }$. Because $M^{2}$ is oriented, we may
choose a global vector field $\eta $ normal along $g$ such that, $%
\left\langle \eta ,w\right\rangle =0$ and $\left\langle \eta ,\eta
\right\rangle =a^{2}-1$. Consider the vector fields 
\begin{equation*}
e_{3}:=\frac{a\eta -w}{a^{2}-1},\text{ \ }e_{4}:=\frac{aw-\eta }{a^{2}-1}%
\text{.}
\end{equation*}
Note that $\left\{ e_{3},e_{4}\right\} $ is a parallel orthonormal frame
field of the normal bundle of $g$ and $e_{4}$ is timelike. Moreover, $%
w=e_{3}+ae_{4}$ and $A_{3}=-aA_{4}$, where $A_{3}$ and $A_{4}$ are the shape
operators of $g$ with respect to the directions $e_{3}$ and $e_{4}$. Since
the second fundamental form of $g$ becomes 
\begin{equation*}
II\left( X,Y\right) =-\left\langle A_{4}X,Y\right\rangle \left(
ae_{3}+e_{4}\right) ,
\end{equation*}
where $X,Y$ are tangent vector fields of $M^{2}$, we deduce that $A_{4}$ is
everywhere nonsingular.$\ $Define now the map $h:M^{2}\rightarrow \mathbb{H}%
^{4}$, $h\left( x\right) :=e_{4}\left( x\right) $, $x\in M^{2}$.\ We claim
that $h$ satisfies all the desired properties. Indeed, for each tangent
vector $X$ of $M^{2}$ we have 
\begin{equation}
dh\left( X\right) =\overset{g}{\nabla }_{X}e_{4}=-dg\left( A_{4}X\right) +%
\overset{g}{D}_{X}e_{4}=-dg\left( A_{4}X\right) \text{.}
\end{equation}
Therefore $h$ is an\ immersion and the metric $\left\langle X,Y\right\rangle
_{h}=\left\langle A_{4}^{2}X,Y\right\rangle $ induced on $M^{2}$ by $h$ is
Riemannian. The normal bundle of $h$ is spanned by $\left\{ \eta
_{1}:=g,\eta _{2}:=e_{3}\right\} $. Denote by $\widetilde{A}_{1}$ and $%
\widetilde{A}_{2}$ the corresponding shape operators of $h$ in the
directions $\eta _{1}$ and $\eta _{2}$. Then the Weingarten formula and $%
\left( 2.6\right) $ yield 
\begin{equation*}
dg\left( X\right) =\overset{h}{\nabla }_{X}\eta _{1}=-dh\left( \widetilde{A}%
_{1}X\right) +\overset{h}{D}_{X}\eta _{1}=dg\left( A_{4}\widetilde{A}%
_{1}X\right) +\overset{h}{D}_{X}\eta _{1}\text{.}
\end{equation*}
Hence 
\begin{equation}
\widetilde{A}_{1}=A_{4}^{-1}\text{ \ and \ }\overset{h}{D}_{X}\eta _{1}=0.
\end{equation}
Moreover, 
\begin{equation*}
d\eta _{2}\left( X\right) =\overset{h}{\nabla }_{X}\eta _{2}=-dh\left( 
\widetilde{A}_{2}X\right) +\overset{h}{D}_{X}\eta _{2}=dg\left( A_{4}%
\widetilde{A}_{2}X\right) +\overset{h}{D}_{X}\eta _{2}\text{.}
\end{equation*}
Since, 
\begin{equation*}
d\eta _{2}\left( X\right) =\overset{g}{\nabla }_{X}e_{3}=-dg\left(
A_{3}X\right) ,
\end{equation*}
we get 
\begin{equation}
\widetilde{A}_{2}=aI\text{ \ and \ }\overset{h}{D}_{X}\eta _{2}=0.
\end{equation}
From $\left( 2.7\right) $ and $\left( 2.8\right) $ we deduce that the vector
field $w=\eta _{2}+ah$ is constant along $h$ and $\left\langle
h,w\right\rangle =-a$. Therefore, $h\left( M^{2}\right) $ is contained in an
umbilical hypersurface $Q^{3}$ of $\mathbb{H}^{4}$, with sectional curvature 
$K_{Q^{3}}=-1+a^{2}=-\left\langle w,w\right\rangle $. Furthermore, $%
h:M^{2}\rightarrow Q^{3}$ is a minimal immersion with normal $\eta _{1}$ and 
$g=\widehat{h}$.

\textbf{Case 2}. Assume now that $w$ is null. According to Proposition \ref
{p2}$\left( iii\right) $ $g$ is superminimal. Because $M^{2}$ is oriented we
may choose a global null vector field $\eta $ in the normal bundle of $g$,
such that $\left\langle \eta ,w\right\rangle =1/2$. Define now the vector
fields 
\begin{equation*}
e_{3}:=\eta +w,\text{ \ }e_{4}:=w-\eta \text{.}
\end{equation*}
Obviously $2w=e_{3}+e_{4}$ and $A_{3}=-A_{4}$, where $A_{3}$ and $A_{4}$ are
the shape operators of $g$ with respect to the directions $e_{3}$ and $e_{4}$%
. Note that $\left\{ e_{3},e_{4}\right\} $ is a parallel orthonormal frame
field of the normal bundle of $g$ and $e_{4}$ is timelike. Since the second
fundamental form of $g$ becomes 
\begin{equation*}
II\left( X,Y\right) =-\left\langle A_{4}X,Y\right\rangle \left(
e_{3}+e_{4}\right) ,
\end{equation*}
where $X,Y$ are tangent vector fields of $M^{2}$, it follows that $A_{4}$ is
everywhere nonsingular. Consider now the map $h:M^{2}\rightarrow \mathbb{H}%
^{4}$, $h\left( x\right) :=e_{4}\left( x\right) $, $x\in M^{2}$. We claim
that $h$ is the required map. Indeed, for each tangent vector $X$ of $M^{2}$
we have 
\begin{equation*}
dh\left( X\right) =-dg\left( A_{4}X\right) .
\end{equation*}
Therefore, $h$ is an immersion and the metric $\left\langle X,Y\right\rangle
_{h}=\left\langle A_{4}^{2}X,Y\right\rangle $ induced by $h$ in $M^{2}$ is
Riemannian. Then the rest of the proof proceeds as in Case 1.
\end{proof}

Let $M^{2}$ be an oriented, 2-dimensional Riemannian manifold and $%
h:M^{2}\rightarrow Q^{3}$ an isometric\ immersion, where $Q^{3}$ is an
umbilical hypersurface of $\mathbb{H}^{4}$. Denote by $\eta $ a unit normal
vector field of $Q^{3}$ in $\mathbb{H}^{4}$. Consider the map $%
F_{h}:M^{2}\times \mathbb{R\rightarrow H}^{4}$, given by 
\begin{equation*}
F_{h}\left( x,t\right) =\cosh th\left( x\right) +\sinh t\eta \circ h\left(
x\right) ,\ \ \left( x,t\right) \in M^{2}\times \mathbb{R},
\end{equation*}
which is called \textit{the suspension }of $h$ in $\mathbb{H}^{4}$.

It is clear that in the case where $h:M^{2}\rightarrow Q^{3}$ is a minimal
isometric immersion without totally geodesic points, then $F_{h}=\Psi _{%
\widehat{h}}\circ T$, where $T$ is the diffeomorphism $T:M^{2}\times \mathbb{%
R\rightarrow }\mathcal{N}^{1}(\widehat{h})$ given by $T\left( x,t\right)
=\left( x,F_{h}\left( x,t\right) \right) $.

In the following proposition we show that there is an abundance of complete
minimal hypersurfaces in $\mathbb{H}^{4}$ with Gauss-Kronecker curvature
identically zero.

\begin{proposition}
\label{p4}Let $h:M^{2}\rightarrow Q^{3}$ be a minimal isometric immersion of
a 2-dimensional, oriented Riemannian manifold $M^{2}$ into an umbilical
hypersurface $Q^{3}$ of $\mathbb{H}^{4}$. Then,

\begin{enumerate}
\item[$\left( i\right) $]  On the open subset of the\ regular points, the
suspension $F_{h}$ of $h$ is a minimal immersion in $\mathbb{H}^{4}$ with
Gauss-Kronecker curvature identically zero.

\item[$\left( ii\right) $]  The metric induced on $M^{2}\times \mathbb{R}$
by $F_{h}$ is complete if and only if $M^{2}$ is\ complete and $Q^{3}$ is a
horosphere or an\ equidistant hypersurface in $\mathbb{H}^{4}$.
\end{enumerate}
\end{proposition}

\begin{proof}
$\left( i\right) $ Denote by $N$ a unit normal vector field along $h$ in $%
Q^{3}$ and by $A_{N}$ the corresponding shape operator of $h$. Let $A_{\eta
}=\alpha I$, $\alpha \in \mathbb{R}$,\ denotes the shape operator of $Q^{3}$
in $\mathbb{H}^{4}$ with respect to the unit normal vector field $\eta $.
Then 
\begin{equation*}
dF_{h}\left( \partial /\partial t\right) =\sinh th+\cosh t\eta \circ h,
\end{equation*}
and for each tangent vector $X$ of $M^{2}$, we have 
\begin{eqnarray*}
dF_{h}\left( X\right) &=&\cosh tdh\left( X\right) +\sinh td\eta \left(
dh\left( X\right) \right) \\
&=&\cosh tdh\left( X\right) -\sinh tA_{\eta }\left( dh\left( X\right) \right)
\\
&=&\left( \cosh t-\alpha \sinh t\right) dh\left( X\right) .
\end{eqnarray*}
Therefore, the point $\left( x,t\right) $ is a regular point of $F_{h}$ if
and only $\cosh t-\alpha \sinh t\neq 0$. The unit vector field $\xi $ given
by $\xi \left( x,t\right) =N\left( x\right) $, $\left( x,t\right) \in
M^{2}\times \mathbb{R}$, is normal along $F_{h}$. Denote by $A_{\xi }$ the
corresponding shape operator. Then $d\xi \left( \frac{\partial }{\partial t}%
\right) =0$, and for each tangent vector $X$ of $M^{2}$, we get 
\begin{eqnarray*}
dF_{h}\left( A_{\xi }X\right) &=&-d\xi \left( X\right) =-dN\left( X\right)
=dh\left( A_{N}X\right) \\
&=&\frac{1}{\cosh t-\alpha \sinh t}dF_{h}\left( A_{N}X\right) .
\end{eqnarray*}
Hence, the principal curvatures of $F_{h}$ are 
\begin{equation}
k_{1}\left( x,t\right) =-k_{3}\left( x,t\right) =\frac{k\left( x\right) }{%
\cosh t-\alpha \sinh t},\ k_{2}\left( x,t\right) =0,
\end{equation}
where $k$ is a principal curvature of $h$.

$\left( ii\right) $ The map $F_{h}$ is an immersion if and only if $\cosh
t-\alpha \sinh t\neq 0$, for each $t\in \mathbb{R}$. This holds if and only
if $\alpha ^{2}\leq 1$. Since the sectional curvature of $Q^{3}$ is $%
K_{Q^{3}}=-1+\alpha ^{2}$, we deduce that the map $F_{h}$ is an immersion if
and only if $Q^{3}$ is a horosphere or an equidistant hypersurface in $%
\mathbb{H}^{4}$. Furthermore, the metric $\left\langle \text{ },\text{ }%
\right\rangle _{F_{h}}$ induced on $M^{2}\times \mathbb{R}$ by $F_{h}$, is
the warped product 
\begin{equation*}
\left\langle \text{ },\text{ }\right\rangle _{F_{h}}=dt^{2}+\left( \cosh
t-\alpha \sinh t\right) ^{2}\left\langle \text{ },\text{ }\right\rangle ,
\end{equation*}
where $\left\langle \text{ },\text{ }\right\rangle $ is the Riemannian
metric of $M^{2}$. Appealing to a result due to Bishop and O'Neill (\cite{BO}%
, Lemma 7.2), $\left\langle \text{ },\text{ }\right\rangle _{F_{h}}$ is
complete if and only if $\left\langle \text{ },\text{ }\right\rangle $ is
complete.
\end{proof}

\begin{remark}
The polar map associated with a non complete stationary surface in $\mathbb{S%
}_{1}^{4}$ may gives rise to a complete minimal hypersurface in $\mathbb{H}%
^{4}$ with Gauss-Kronecker curvature identically zero and nowhere vanishing
second fundamental form. In fact, consider a $2$-dimensional, oriented,
complete Riemannian manifold $M^{2}$ and suppose that $h:M^{2}\rightarrow
Q^{3}$ is a minimal isometric immersion without totally geodesic points,
where $Q^{3}$ is a horosphere or an equidistant hypersurface of $\mathbb{H}%
^{4}$. The metric $\left\langle \text{ , }\right\rangle _{\widehat{h}}$
induced by $\widehat{h}$ is not complete. Indeed, if $\left\langle \text{ , }%
\right\rangle _{\widehat{h}}$ where complete, then by Myers' Theorem and the
fact that its Gaussian curvature $K_{\widehat{h}}$ satisfies 
\begin{equation*}
K_{\widehat{h}}=1-\frac{K_{Q^{3}}}{K_{Q^{3}}-K}\geq 1,
\end{equation*}
$M^{2}$ would be compact. A contradiction, since there are no compact
minimal surfaces in simply connected space forms of non positive sectional
curvature. Moreover according to Proposition \ref{p4}, the metric induced on 
$\mathcal{N}^{1}(\widehat{h})$ by $\Psi _{\widehat{h}}$ is\ complete.
\end{remark}

\begin{remark}
There are numerous examples of complete minimal hypersurfaces in $\mathbb{H}%
^{4}$ with Gauss-Kronecker curvature identically zero and unbounded scalar
curvature. Indeed, suppose that $Q^{3}$ is a horosphere and $%
h:M^{2}\rightarrow Q^{3}$ is a complete minimal immersion. Then the
suspension of $h$ is a complete hypersurface in $\mathbb{H}^{4}$. According
to $\left( 2.9\right) $ its principal curvatures are 
\begin{equation*}
k_{1}\left( x,t\right) =-k_{3}\left( x,t\right) =\frac{k\left( x\right) }{%
\cosh t-\sinh t}\text{, }k_{2}\left( x,t\right) =0\text{.}
\end{equation*}
Because $\lim\limits_{t\rightarrow \infty }k_{1}\left( x,t\right) =\infty $,
it follows that the scalar curvature of the suspension must be unbounded.
There are also plenty of complete minimal hypersurfaces in $\mathbb{H}^{4}$
with Gauss-Kronecker curvature zero and bounded scalar curvature. Indeed,
there exist complete minimal surfaces in $\mathbb{H}^{3}$ with Gaussian
curvature bounded from below (see \cite{CD}). The suspension of such
surfaces are minimal hypersurfaces in $\mathbb{H}^{4}$ with Gauss-Kronecker
curvature identically zero and bounded scalar curvature.
\end{remark}

\section{Local theory of minimal hypersurfaces in $\mathbb{H}^{4}$ with zero
Gauss-Kronecker curvature}

Let $M^{3}$ be a 3-dimensional, oriented Riemannian manifold and $%
f:M^{3}\rightarrow \mathbb{H}^{4}$ an isometric minimal immersion. Denote by 
$\xi $ a unit normal vector field along $f$ with corresponding shape
operator $A$ and principal curvatures $k_{1}\geq k_{2}\geq k_{3}$. The 
\textit{Gauss-Kronecker curvature} $\mathcal{K}$ of $f$ and the \textit{%
scalar curvature} $\tau $ of $M^{3}$ are given by 
\begin{equation*}
\mathcal{K}=k_{1}k_{2}k_{3},\ \tau =-6-\left(
k_{1}^{2}+k_{2}^{2}+k_{3}^{2}\right) .
\end{equation*}
Assume now that the second fundamental form of $f$ is nowhere zero and that
the Gauss-Kronecker curvature is identically zero. Then the principal
curvatures are $k_{1}=\lambda $, $k_{2}=0$, $k_{3}=-\lambda $, where $%
\lambda $ is a smooth positive function on $M^{3}$. We can choose locally an
orthonormal frame field $\left\{ e_{1},e_{2},e_{3}\right\} $ of principal
directions corresponding to $\lambda ,0,-\lambda $. Let $\left\{ \omega
_{1},\omega _{2},\omega _{3}\right\} $ and $\left\{ \omega _{ij}\right\} $, $%
i,j\in \{1,2,3\}$, be the dual frame and the\ connection forms. Hereafter we
make the following convection on the ranges of indices 
\begin{equation*}
1\leq i,j,k,\ldots \leq 3,
\end{equation*}
and adopt the method of moving frames. The structure equations are 
\begin{eqnarray*}
d\omega _{i} &=&\sum_{j}\omega _{ij}\wedge \omega _{j},\quad \omega
_{ij}+\omega _{ji}=0, \\
d\omega _{ij} &=&\sum_{l}\omega _{il}\wedge \omega _{lj}-\left(
k_{i}k_{j}-1\right) \omega _{i}\wedge \omega _{j}.
\end{eqnarray*}
Consider the functions 
\begin{equation*}
u:=\omega _{12}\left( e_{3}\right) ,\quad v:=e_{2}\left( \log \lambda
\right) ,
\end{equation*}
which will play a crucial role in the sequel. From the structural equations,
and the Codazzi equations, 
\begin{equation*}
\begin{array}{l}
e_{i}\left( k_{j}\right) =\left( k_{i}-k_{j}\right) \omega _{ij}\left(
e_{j}\right) ,\ i\neq j,\medskip \\ 
\left( k_{1}-k_{2}\right) \omega _{12}\left( e_{3}\right) =\left(
k_{2}-k_{3}\right) \omega _{23}\left( e_{1}\right) =\left(
k_{1}-k_{3}\right) \omega _{13}\left( e_{2}\right) ,
\end{array}
\end{equation*}
we easily get 
\begin{equation}
\begin{array}{lll}
\omega _{12}\left( e_{1}\right) =v, & \omega _{13}\left( e_{1}\right) =\frac{%
1}{2}e_{3}\left( \log \lambda \right) , & \omega _{23}\left( e_{1}\right)
=u,\medskip \\ 
\omega _{12}\left( e_{2}\right) =0, & \omega _{13}\left( e_{2}\right) =\frac{%
1}{2}u, & \omega _{23}\left( e_{2}\right) =0,\medskip \\ 
\omega _{12}\left( e_{3}\right) =u, & \omega _{13}\left( e_{3}\right) =-%
\frac{1}{2}e_{1}\left( \log \lambda \right) , & \omega _{23}\left(
e_{3}\right) =-v
\end{array}
\end{equation}
and 
\begin{equation}
\begin{array}{l}
e_{2}\left( v\right) =v^{2}-u^{2}-1,\ e_{1}\left( u\right) =e_{3}\left(
v\right) ,\medskip \\ 
e_{2}\left( u\right) =2uv,\ e_{3}\left( u\right) =-e_{1}\left( v\right) .
\end{array}
\end{equation}
Furthermore, the\ above equations yield 
\begin{equation}
\begin{array}{c}
\lbrack e_{1},e_{2}]=-ve_{1}+\frac{1}{2}ue_{3},\ [e_{2},e_{3}]=\frac{1}{2}%
ue_{1}+ve_{3},\medskip \\ 
\lbrack e_{1},e_{3}]=-\frac{1}{2}e_{3}\left( \log \lambda \right)
e_{1}-2ue_{2}+\frac{1}{2}e_{1}\left( \log \lambda \right) e_{3}.
\end{array}
\end{equation}

\begin{lemma}
\label{l1}The functions $u$ and $v$ are harmonic.
\end{lemma}

\begin{proof}
Using $\left( 3.1\right) $, from the definition of the Laplacian we have 
\begin{eqnarray*}
\Delta v &=&e_{1}e_{1}\left( v\right) +e_{2}e_{2}\left( v\right)
+e_{3}e_{3}\left( v\right) -\left( \omega _{31}\left( e_{3}\right) +\omega
_{21}\left( e_{2}\right) \right) e_{1}\left( v\right) \\
&&-\left( \omega _{12}\left( e_{1}\right) +\omega _{32}\left( e_{3}\right)
\right) e_{2}\left( v\right) -\left( \omega _{13}\left( e_{1}\right) +\omega
_{23}\left( e_{2}\right) \right) e_{3}\left( v\right) \\
&=&e_{1}e_{1}\left( v\right) +e_{2}e_{2}\left( v\right) +e_{3}e_{3}\left(
v\right) -\frac{1}{2}e_{1}\left( \log \lambda \right) e_{1}\left( v\right)
-2ve_{2}\left( v\right) \\
&&-\frac{1}{2}e_{3}\left( \log \lambda \right) e_{3}\left( v\right) \text{.}
\end{eqnarray*}
In view of $\left( 3.2\right) $, we get 
\begin{eqnarray*}
e_{1}e_{1}\left( v\right) &=&-e_{1}e_{3}\left( u\right) ,\ e_{3}e_{3}\left(
v\right) =e_{3}e_{1}\left( u\right) , \\
e_{2}e_{2}\left( v\right) &=&2ve_{2}\left( v\right) -2ue_{2}\left( u\right) 
\text{.}
\end{eqnarray*}
On account of $\left( 3.2\right) $, $\left( 3.3\right) $ and the previous
relations, we obtain 
\begin{eqnarray*}
\Delta v &=&-e_{1}e_{3}\left( u\right) +e_{3}e_{1}\left( u\right)
+2ve_{2}\left( v\right) -2ue_{2}\left( u\right) \\
&&-\frac{1}{2}e_{1}\left( \log \lambda \right) e_{1}\left( v\right)
-2ve_{2}\left( v\right) -\frac{1}{2}e_{3}\left( \log \lambda \right)
e_{3}\left( v\right) \\
&=&\frac{1}{2}e_{3}\left( \log \lambda \right) e_{1}\left( u\right)
+2ue_{2}\left( u\right) -\frac{1}{2}e_{1}\left( \log \lambda \right)
e_{3}\left( u\right) \\
&&-2ue_{2}\left( u\right) -\frac{1}{2}e_{1}\left( \log \lambda \right)
e_{1}\left( v\right) -\frac{1}{2}e_{3}\left( \log \lambda \right)
e_{3}\left( v\right) \\
&=&0\text{.}
\end{eqnarray*}
In a similar way, we verify that $\Delta u=0$.
\end{proof}

\begin{lemma}
\label{l2}Let $\gamma :I\subset \mathbb{R\rightarrow }M^{3}$ be an integral
curve of $e_{2}$ emanating from $x\in M^{3}$. Then $\gamma $ is a\ geodesic
of $M^{3}$ and $f\circ \gamma $ is a geodesic of $\mathbb{H}^{4}$. Moreover, 
\begin{equation*}
\frac{1}{\lambda ^{2}\circ \gamma \left( t\right) }=\frac{1}{2}\left(
a\left( x\right) e^{2t}+b\left( x\right) +d\left( x\right) e^{-2t}\right)
\end{equation*}
and 
\begin{equation*}
v\circ \gamma \left( t\right) =-\dfrac{a\left( x\right) e^{2t}-d\left(
x\right) e^{-2t}}{a\left( x\right) e^{2t}+b\left( x\right) +d\left( x\right)
e^{-2t}},
\end{equation*}
where $a\left( x\right) ,$ $b\left( x\right) ,$ $d\left( x\right) $ are real
constants depending only on $x$ and $t\in I$.
\end{lemma}

\begin{proof}
By making use of $\left( 3.1\right) $ we, immediately, obtain $\nabla
_{e_{2}}e_{2}=0$. Thus, $\gamma $ is a geodesic of $M^{3}$ and the Gauss
formula implies that $f\circ \gamma $ is a geodesic of $\mathbb{H}^{4}$. By
virtue of $\left( 3.2\right) $, we easily get 
\begin{equation*}
e_{2}e_{2}e_{2}\left( \frac{1}{\lambda ^{2}}\right) =4e_{2}\left( \frac{1}{%
\lambda ^{2}}\right) .
\end{equation*}
Restricting the last equation along $\gamma $ and integrating, we deduce
that 
\begin{equation*}
\frac{1}{\lambda ^{2}\circ \gamma \left( t\right) }=\frac{1}{2}\left(
a\left( x\right) e^{2t}+b\left( x\right) +d\left( x\right) e^{-2t}\right) ,
\end{equation*}
where $a\left( x\right) ,$ $b\left( x\right) ,$ $d\left( x\right) $ are real
constants. Differentiating, we obtain 
\begin{equation*}
v\circ \gamma \left( t\right) =\frac{d}{dt}\left( \log \lambda \circ \gamma
\right) \left( t\right) =-\dfrac{a\left( x\right) e^{2t}-d\left( x\right)
e^{-2t}}{a\left( x\right) e^{2t}+b\left( x\right) +d\left( x\right) e^{-2t}},
\end{equation*}
and the proof is finished.
\end{proof}

We are now ready to give the local classification of minimal hypersurfaces
in $\mathbb{H}^{4}$ with Gauss-Kronecker curvature zero and nowhere
vanishing second fundamental form, which is in fact the converse of
Proposition \ref{p1}$\left( ii\right) $.

\begin{proposition}
\label{p5}Let $M^{3}$ be a 3-dimensional, oriented, Riemannian manifold and $%
f:M^{3}\rightarrow \mathbb{H}^{4}$ a minimal isometric immersion with unit
normal vector field $\xi $, Gauss-Kronecker curvature identically zero and
nowhere vanishing second fundamental form. Each point\textit{\ }$x_{0}\in
M^{3}$ has a neighborhood $U$ such that the quotient space $V$ of leaves of
the nullity distribution on $U$ is a 2-dimensional differentiable manifold
with quotient projection $\pi :U\rightarrow V$ and

\begin{enumerate}
\item[$\left( i\right) $]  there exists a spacelike stationary immersion $%
g:V\rightarrow \mathbb{S}_{1}^{4}$ and an isometry $T:U\rightarrow \mathcal{N%
}^{1}\left( g\right) $ such that $g\circ \pi =\xi $ and $f=\Psi _{g}\circ T$
on $U$,

\item[$\left( ii\right) $]  the Gaussian curvature $K$ of the metric induced
by $g$ on $V$ and the normal curvature $K^{\bot }$ satisfy 
\begin{equation*}
K\circ \pi =1+\frac{1-u^{2}-v^{2}}{\lambda ^{2}},\ \ K^{\bot }\circ \pi =-%
\frac{2u}{\lambda ^{2}}
\end{equation*}
on $U$.
\end{enumerate}
\end{proposition}

\begin{proof}
Consider a coordinate system $\left( x_{1},x_{2},x_{3}\right) $ on $U\subset
M^{3}$, around $x_{0}$, such that $\frac{\partial }{\partial x_{2}}=e_{2}$.
Denote by $V$ the quotient space of leaves of the nullity distribution on $U$
and\ by $\pi :U\rightarrow V$ the quotient projection. It is well known that 
$V$ can be equipped with a structure of a 2-dimensional differentiable
manifold which makes $\pi $ a submersion.

Our assumptions ensure that the unit normal vector field $\xi $ remains
constant along each leaf of the nullity distribution and so we may define a
smooth map $g:V\rightarrow \mathbb{S}_{1}^{4}$ so that $g\circ \pi =\xi $.
We claim that $g$ is a spacelike stationary immersion. Indeed, consider a
smooth transversal $S$ to the leaves of the nullity distribution, through a
point $x\in U$ such that the frame $\left\{ E_{1}:=e_{1}\left( x\right)
,E_{3}:=e_{3}\left( x\right) \right\} $ spans $T_{x}S$. Because $\pi $ is a
submersion,$\ \left\{ d\pi \left( E_{1}\right) ,d\pi \left( E_{3}\right)
\right\} $ constitute a base of $T_{\pi \left( x\right) }V$. Note that 
\begin{equation*}
dg\left( d\pi \left( E_{1}\right) \right) =-\lambda \left( x\right) df\left(
E_{1}\right) \text{ and }dg\left( d\pi \left( E_{3}\right) \right) =\lambda
\left( x\right) df\left( E_{3}\right) .
\end{equation*}
Thus $g$ is a spacelike immersion and $\left\{ X_{1}:=\frac{1}{\lambda
\left( x\right) }d\pi \left( E_{1}\right) ,X_{3}:=\frac{1}{\lambda \left(
x\right) }d\pi \left( E_{3}\right) \right\} $ is an orthonormal base at $\pi
\left( x\right) $ with respect to the metric induced by $g$. Let $\left\{
\eta _{3},\eta _{4}\right\} $ be an orthonormal frame field in the normal
bundle of $g$ such that $\eta _{3}\circ \pi =df\left( e_{2}\right) $ and $%
\eta _{4}\circ \pi =f$ on $S$. Then bearing in mind the Gauss formula and $%
\left( 3.1\right) $ we obtain\ 
\begin{eqnarray*}
d\eta _{3}\left( X_{1}\right) &=&-\frac{v\left( x\right) }{\lambda \left(
x\right) }df\left( E_{1}\right) +\frac{u\left( x\right) }{\lambda \left(
x\right) }df\left( E_{3}\right) \medskip , \\
d\eta _{3}\left( X_{3}\right) &=&-\frac{u\left( x\right) }{\lambda \left(
x\right) }df\left( E_{1}\right) -\frac{v\left( x\right) }{\lambda \left(
x\right) }df\left( E_{3}\right) \medskip , \\
d\eta _{4}\left( X_{1}\right) &=&\frac{1}{\lambda \left( x\right) }df\left(
E_{1}\right) \text{ \ and \ }d\eta _{4}\left( X_{3}\right) =\frac{1}{\lambda
\left( x\right) }df\left( E_{3}\right) .
\end{eqnarray*}
Denote by $A_{3}$, $A_{4}$ the shape operators of $g$ at $\pi \left(
x\right) $ corresponding to the normal directions $\eta _{3}$ and $\eta _{4}$%
. Taking into account the above relations, from Weingarten formulas it
follows that at $\pi \left( x\right) $ we have 
\begin{equation}
A_{3}\sim \frac{1}{\lambda \left( x\right) }\left( 
\begin{array}{ll}
-v\left( x\right) & -u\left( x\right) \\ 
-u\left( x\right) & \ \ v\left( x\right)
\end{array}
\right) ,\quad A_{4}\sim \frac{1}{\lambda \left( x\right) }\left( 
\begin{array}{ll}
1 & \ \ 0 \\ 
0 & -1
\end{array}
\right) ,
\end{equation}
with respect to the orthonormal base $\left\{ X_{1},X_{3}\right\} $. So the
immersion $g:V\rightarrow \mathbb{S}_{1}^{4}$ is a stationary immersion.
Moreover we have $f=\Psi _{g}\circ T$ on $U$, where the map $T:U\rightarrow 
\mathcal{N}^{1}\left( g\right) $ is defined by $T\left( x\right) =\left( \pi
\left( x\right) ,f\left( x\right) \right) $, $x\in U$. By restricting $U$,
if necessary, $T$ is an isometry because $\mathcal{N}^{1}\left( g\right) $
is equipped with the metric induced by $\Psi _{g}$.

Part $\left( ii\right) $ follows immediately from $\left( 3.4\right) $.
\end{proof}

\section{Complete minimal hypersurfaces in $\mathbb{H}^{4}$ with vanishing
Gauss-Kronecker curvature}

The purpose of this section is to classify complete minimal hypersurfaces in 
$\mathbb{H}^{4}$ with Gauss-Kronecker curvature identically zero and nowhere
zero second fundamental form, under the assumption that the scalar curvature
is bounded from below. More precisely, we shall prove the following

\begin{Theor.}
Let $M^{3}$ be a 3-dimensional, oriented,\ complete Riemannian manifold
whose scalar curvature is bounded from below and $f:M^{3}\rightarrow \mathbb{%
H}^{4}$ a minimal isometric immersion with Gauss-Kronecker curvature
identically zero and nowhere zero second fundamental form. Then there exist
a minimal isometric immersion $h:M^{2}\rightarrow Q^{3}$, without totally
geodesic points, of a complete 2-dimensional oriented Riemannian manifold $%
M^{2}$ into an equidistant hypersurface $Q^{3}$ of $\mathbb{H}^{4}$ and a
local isometry $T:M^{3}\rightarrow \mathcal{N}^{1}(\widehat{h})$ such that $%
f=\Psi _{\widehat{h}}\circ T$.
\end{Theor.}

The proof of our result relies heavily on the well known Generalized Maximum
Principle due to Omori and\ Yau (\cite{O},\cite{Y}):

\begin{GMP}
Let $M$\ be a complete Riemannian manifold whose Ricci curvature is bounded
from below. If $\varphi $\ is a $C^{2}$-function on $M$ bounded from above,
then there exists a sequence $\left\{ x_{n}\right\} $\ of points of $M$\
such that 
\begin{equation*}
\lim \varphi \left( x_{n}\right) =\sup \varphi ,\text{ }\left| \nabla
\varphi \right| \left( x_{n}\right) \leq \frac{1}{n}\ \text{and }\Delta
\varphi \left( x_{n}\right) \leq \frac{1}{n},
\end{equation*}
for each $n\in \mathbb{N}$, where $\nabla $, $\Delta $ stand for the
gradient and Laplacian operator.
\end{GMP}

The following lemma, is essentially a consequence of a result proved by
Cheng and Yau (\cite{CY}, Theorem 8). For reader's convenience we shall give
here a short proof.

\begin{lemma}
\label{l3}Let $M$ be a complete Riemannian manifold whose Ricci curvature is
bounded from below, and $\varphi $ a $C^{2}$-solution of the differential
inequality 
\begin{equation*}
\Delta \varphi \geq 2\varphi ^{2}.
\end{equation*}
Then $\varphi $ is bounded from above and $\sup \varphi =0$.
\end{lemma}

\begin{proof}
We suppose in the contrary that $\sup \varphi =\infty $. Then there exists a
point $x_{0}\in M$ such that $\varphi \left( x_{0}\right) \geq 2$. Consider
a $C^{2}$-positive increasing function $F:\mathbb{R\rightarrow R}$ which for 
$t\geq 2$ is given by $F\left( t\right) =2\left( 1-t^{-1/2}\right) $. The
function $\Phi =F\circ \varphi $ is bounded from above, since $\Phi \leq 2$.
Appealing to the Generalized Maximum Principle, we deduce that there exists
a sequence $\left\{ x_{n}\right\} $ such that 
\begin{equation*}
\lim \Phi \left( x_{n}\right) =\sup \Phi ,\text{\ \ }\left| \nabla \Phi
\right| \left( x_{n}\right) \leq \frac{1}{n}\ \text{and}\ \Delta \Phi \left(
x_{n}\right) \leq \frac{1}{n},
\end{equation*}
for each $n\in \mathbb{N}$. For $n$ large enough we have $\varphi \left(
x_{n}\right) \geq 2$. Hence, estimating at $x_{n}$ we get 
\begin{equation}
\left| \nabla \Phi \right| \left( x_{n}\right) =F^{\prime }\left( \varphi
\left( x_{n}\right) \right) \left| \nabla \varphi \right| \left(
x_{n}\right) =\varphi ^{-3/2}\left( x_{n}\right) \left| \nabla \varphi
\right| \left( x_{n}\right) \leq \frac{1}{n}
\end{equation}
and 
\begin{equation}
\Delta \Phi \left( x_{n}\right) =\varphi ^{-3/2}\left( x_{n}\right) \Delta
\varphi \left( x_{n}\right) -\frac{3}{2}\varphi ^{-5/2}\left( x_{n}\right)
\left| \nabla \varphi \right| ^{2}\left( x_{n}\right) \leq \frac{1}{n}.
\end{equation}
Combining $\left( 4.2\right) $ with $\left( 4.1\right) $, and bearing in
mind that $\Delta \varphi \geq 2\varphi ^{2}$, we obtain 
\begin{equation*}
2-\frac{3}{2}\varphi ^{-3}\left( x_{n}\right) \left| \nabla \varphi \right|
^{2}\left( x_{n}\right) \leq \frac{1}{n}\varphi ^{-1/2}\left( x_{n}\right) 
\text{.}
\end{equation*}
Letting $n\rightarrow \infty $, we get a contradiction. Therefore $\varphi $
must be bounded from above. Appealing again to the Generalized Maximum
Principle, and\ bearing in mind that $\Delta \varphi \geq 2\varphi ^{2}$, we
infer that $\sup \varphi =0$.
\end{proof}

\begin{proof}[Proof of the Theorem]
Let $A$ be the shape operator associated with a unit normal $\xi $. Then the
principal curvatures of $f$ are $k_{1}=\lambda $, $k_{2}=0$, $k_{3}=-\lambda 
$, where $\lambda $ is a smooth positive function on $M^{3}$. It is well
known that the nullity distribution $\Delta =\ker A$ is smooth. We
distinguish two cases.

\textbf{Case 1.} We assume that there exists a global unit\ section $e_{2}$
of$\ \Delta $. Then the function $v=e_{2}\left( \log \lambda \right) $ is
globally defined and smooth. Around each point $x\in M^{3}$ we may choose a
neighborhood $U_{x}$ of $x$, vector fields $e_{1}$, $e_{3}$ such that the
orthonormal frame field $\left\{ e_{1},e_{2},e_{3}\right\} $ gives the right
orientation of $M^{3}$ and $Ae_{1}=\lambda e_{1}$, $Ae_{3}=-\lambda e_{3}$
on $U_{x}$. If\ \ for another point $\overline{x}\in M^{3}$ with
corresponding neighborhood $U_{\overline{x}}$ and orthonormal frame field $%
\left\{ \overline{e}_{1},e_{2},\overline{e}_{3}\right\} $ chosen\ as before,
we have $U_{x}\cap U_{\overline{x}}\neq \emptyset $, then either $\overline{e%
}_{1}=e_{1}$ and $\overline{e}_{3}=e_{3}$ or $\overline{e}_{1}=-e_{1}$ and $%
\overline{e}_{3}=-e_{3}$ on $U_{x}\cap U_{\overline{x}}$. Thus $\left\langle
\nabla _{e_{3}}e_{1},e_{2}\right\rangle =\left\langle \nabla _{\overline{e}%
_{3}}\overline{e}_{1},e_{2}\right\rangle $ on $U_{x}\cap U_{\overline{x}}$
and so the local function $u$ introduced in Section 3 can be extended to a
smooth global one.

Our assumptions imply that the Ricci curvature of $M^{3}$ is bounded from
below. Making use of $\left( 3.2\right) $ and the harmonicity of $u$ and $v$
(Lemma \ref{l1}), we obtain 
\begin{eqnarray*}
\frac{1}{2}\Delta \left( u^{2}+v^{2}-1\right) &=&\left| \nabla u\right|
^{2}+\left| \nabla v\right| ^{2} \\
&\geq &\left( e_{2}\left( u\right) \right) ^{2}+\left( e_{2}\left( v\right)
\right) ^{2} \\
&=&4u^{2}v^{2}+\left( v^{2}-u^{2}-1\right) ^{2} \\
&\geq &\left( u^{2}+v^{2}-1\right) ^{2}\text{.}
\end{eqnarray*}
Then, by virtue of Lemma \ref{l3}, we have $\sup \left( u^{2}+v^{2}-1\right)
=0$, which implies $u^{2}+v^{2}\leq 1$.

\textit{Claim}: $u\equiv 0$. At first we shall prove that $v^{2}<1$. Arguing
indirectly, we assume that there exists a point $x_{0}\in M^{3}$ such that $%
\left| v\left( x_{0}\right) \right| =1$. The harmonicity of $v$ and the
maximum principle imply either $v\equiv 1$ or $v\equiv -1$. Then Lemma \ref
{l2} yields $a\left( x_{0}\right) =b\left( x_{0}\right) =0$ or $b\left(
x_{0}\right) =d\left( x_{0}\right) =0$, respectively, and thus $\lambda
^{2}\circ \gamma \left( t\right) $, $t\in \mathbb{R}$, is unbounded, where $%
\gamma $ is the integral curve of $e_{2}$ emanating from the point $x_{0}$.
This contradicts our assumption on the scalar curvature. So $v^{2}<1$. It is
obvious from Lemma \ref{l2}, that on each integral curve of $e_{2}$, the
function $v$ changes sign only once.

Consider, now, the set $v^{-1}\left( 0\right) $. From $\left( 3.2\right) $
we have $e_{2}\left( v\right) =v^{2}-u^{2}-1<0$. Hence $0$ is a regular
value of $v$ and thus $v^{-1}\left( 0\right) $ is an oriented and\ connected
2-dimensional submanifold of $M^{3}$. The map $\rho :v^{-1}\left( 0\right)
\times \mathbb{R\rightarrow }M^{3}$ defined by $\rho \left( x,t\right)
:=\exp _{x}\left( te_{2}\left( x\right) \right) $, where $\exp _{x}$ denotes
the exponential map of $M^{3}$ based on the point $x\in v^{-1}\left(
0\right) $,\ is a\ diffeomorphism. Appealing to Lemma \ref{l2}, we have 
\begin{equation*}
v\circ \rho \left( x,t\right) =-\frac{a\left( x\right) e^{2t}-d\left(
x\right) e^{-2t}}{a\left( x\right) e^{2t}+b\left( x\right) +d\left( x\right)
e^{-2t}},
\end{equation*}
where $a\left( x\right) ,b\left( x\right) ,d\left( x\right) $ are\ smooth
functions on $v^{-1}\left( 0\right) $. Since $v\circ \rho \left( x,0\right)
=0$, we obtain $a\left( x\right) =d\left( x\right) $ for each $x\in
v^{-1}\left( 0\right) $. Hence, 
\begin{equation}
\frac{1}{\lambda ^{2}\circ \rho \left( x,t\right) }=a\left( x\right) \cosh
2t+\frac{b\left( x\right) }{2}
\end{equation}
and 
\begin{equation}
v\circ \rho \left( x,t\right) =-\frac{2a\left( x\right) \sinh 2t}{2a\left(
x\right) \cosh 2t+b\left( x\right) }=-a\left( x\right) \sinh 2t\ \lambda
^{2}\circ \rho \left( x,t\right) .
\end{equation}
From $\left( 4.3\right) $, $\left( 4.4\right) $ and in view of $e_{2}\left(
v\right) =v^{2}-u^{2}-1<0$, we deduce that $a\left( x\right) >0$ for each $%
x\in v^{-1}\left( 0\right) $. Consider now the function $\phi :v^{-1}\left(
0\right) \times \mathbb{R\rightarrow R}$, $\phi \left( x,t\right) =\tanh t$.
Since $d\rho \left( \frac{\partial }{\partial t}\right) =e_{2}$, we have 
\begin{equation}
e_{2}\left( \phi \circ \rho ^{-1}\right) =1-\phi ^{2}\circ \rho ^{-1}.
\end{equation}
Differentiating $\left( 4.4\right) $ with respect to $\frac{\partial }{%
\partial t}$ and making use of $\left( 3.2\right) $ and $\left( 4.4\right) $
we obtain 
\begin{equation}
\frac{\phi \circ \rho ^{-1}}{1+\phi ^{2}\circ \rho ^{-1}}=\frac{-v}{%
1+u^{2}+v^{2}}.
\end{equation}
Obviously we have $v\left( \phi \circ \rho ^{-1}\right) \leq 0$. The
function $G:=u^{2}+\left( v+\phi \circ \rho ^{-1}\right) ^{2}$ is smooth and
bounded from above. Appealing to the Generalized Maximum Principle, there
exists a sequence $\left\{ x_{n}\right\} $ of points in $M^{3}$ such that 
\begin{equation*}
\lim G\left( x_{n}\right) =\sup G,\text{ }\left| \nabla G\right| \left(
x_{n}\right) \leq \frac{1}{n}\text{ and\ }\Delta G\left( x_{n}\right) \leq 
\frac{1}{n},
\end{equation*}
for each $n\in \mathbb{N}$. Because the functions $u$, $v$ and$\ \phi \circ
\rho ^{-1}$ are bounded, without loss of generality, we may assume that 
\begin{equation*}
\lim u\left( x_{n}\right) =u_{0},\text{ }\lim v\left( x_{n}\right) =v_{0}%
\text{ and }\lim \phi \circ \rho ^{-1}\left( x_{n}\right) =\phi _{0},
\end{equation*}
where $u_{0}$, $v_{0}$ and $\phi _{0}$ are real numbers. Using the equations 
$\left( 3.2\right) $, $\left( 4.5\right) $ and the harmonicity of $u$ and $v$%
, we readily see that 
\begin{eqnarray}
\frac{1}{2}e_{2}\left( G\right)  &=&ue_{2}\left( u\right) +\left( v+\phi
\circ \rho ^{-1}\right) \left( e_{2}\left( v\right) +e_{2}\left( \phi \circ
\rho ^{-1}\right) \right)  \\
&=&2u^{2}v+\left( v+\phi \circ \rho ^{-1}\right) \left( v^{2}-u^{2}-\phi
^{2}\circ \rho ^{-1}\right)   \notag \\
&=&\left( v-\phi \circ \rho ^{-1}\right) G  \notag
\end{eqnarray}
and 
\begin{eqnarray}
\frac{1}{2}\Delta G &=&\left| \nabla u\right| ^{2}+\left( v+\phi \circ \rho
^{-1}\right) \Delta \left( \phi \circ \rho ^{-1}\right) +\left| \nabla
\left( v+\phi \circ \rho ^{-1}\right) \right| ^{2} \\
&\geq &4u^{2}v^{2}+\left( v+\phi \circ \rho ^{-1}\right) \Delta \left( \phi
\circ \rho ^{-1}\right)   \notag \\
&&+\left( e_{2}\left( v\right) +e_{2}\left( \phi \circ \rho ^{-1}\right)
\right) ^{2}  \notag \\
&=&4u^{2}v^{2}+\left( v+\phi \circ \rho ^{-1}\right) \Delta \left( \phi
\circ \rho ^{-1}\right)   \notag \\
&&+\left( v^{2}-u^{2}-\phi ^{2}\circ \rho ^{-1}\right) ^{2}.  \notag
\end{eqnarray}
Estimating at $x_{n}$ and letting $n\rightarrow \infty $,$\ $the equation$\
\left( 4.7\right) $ yields 
\begin{equation*}
\left( v_{0}-\phi _{0}\right) \sup G=0.
\end{equation*}
If $v_{0}\neq \phi _{0}$ we obtain $\sup G=0$, which proves our claim.

Suppose now that $v_{0}=\phi _{0}$. Then, because of $v\left( \phi \circ
\rho ^{-1}\right) \leq 0$, we get $v_{0}=\phi _{0}=0$. Making use of $\left(
4.6\right) $ and of\ the harmonicity of the functions$\ u$ and $v$, a
straightforward computation\ shows that 
\begin{equation}
\dfrac{1-\phi ^{2}\circ \rho ^{-1}}{\left( 1+\phi ^{2}\circ \rho
^{-1}\right) ^{2}}\nabla \left( \phi \circ \rho ^{-1}\right) =\dfrac{2uv}{%
\left( 1+u^{2}+v^{2}\right) ^{2}}\nabla u-\frac{1+u^{2}-v^{2}}{\left(
1+u^{2}+v^{2}\right) ^{2}}\nabla v
\end{equation}
and 
\begin{equation}
\begin{array}{l}
\dfrac{1-\phi ^{2}\circ \rho ^{-1}}{\left( 1+\phi ^{2}\circ \rho
^{-1}\right) ^{2}}\Delta \left( \phi \circ \rho ^{-1}\right) =\dfrac{%
2v\left( 1-3u^{2}+v^{2}\right) }{\left( 1+u^{2}+v^{2}\right) ^{3}}\left|
\nabla u\right| ^{2} \\ 
\qquad \qquad \qquad \quad +\dfrac{4u\left( 1+u^{2}-3v^{2}\right) }{\left(
1+u^{2}+v^{2}\right) ^{3}}\left\langle \nabla u,\nabla v\right\rangle +%
\dfrac{2v\left( 3+3u^{2}-v^{2}\right) }{\left( 1+u^{2}+v^{2}\right) ^{3}}%
\left| \nabla v\right| ^{2} \\ 
\qquad \qquad \qquad \quad +2\left( \phi \circ \rho ^{-1}\right) \dfrac{%
3-\phi ^{2}\circ \rho ^{-1}}{\left( 1+\phi ^{2}\circ \rho ^{-1}\right) ^{3}}%
\left| \nabla \left( \phi \circ \rho ^{-1}\right) \right| ^{2}.
\end{array}
\end{equation}
Since $u$, $v$ are bounded harmonic functions and $M^{3}$ has Ricci
curvature bounded from below, by a result due to Yau (\cite{Y}, Theorem $%
3^{\prime \prime }$), it follows that the functions $\left| \nabla u\right|
^{2}$ and $\left| \nabla v\right| ^{2}$ are also bounded. Hence, from $%
\left( 4.9\right) $ and $\left( 4.10\right) $ we deduce that the sequence $%
\left\{ \Delta \left( \phi \circ \rho ^{-1}\right) \left( x_{n}\right)
\right\} $ is bounded. Estimating at $x_{n}$ and passing to the limit, from$%
\ \left( 4.8\right) $ we get $u_{0}=0$. So $\sup G=0$, because of $%
v_{0}=\phi _{0}=0$. Thus $G\equiv 0$ and consequently $u\equiv 0$, which
completes the proof of our claim.

It can be easily seen that the quotient space $M^{2}$ of leaves of the
nullity distribution can be identified with the manifold $v^{-1}\left(
0\right) $ via the diffeomorphism $\rho $. Hence, $M^{2}$ inherits in a
natural way the structure of a 2-dimensional manifold that makes the
quotient projection $\pi :M^{3}\rightarrow M^{2}$ a submersion. Thus
appealing to Proposition \ref{p5}, there exists a spacelike stationary
immersion $g:M^{2}\rightarrow \mathbb{S}_{1}^{4}$ and an isometry $%
T:M^{3}\rightarrow \mathcal{N}^{1}\left( g\right) $, defined by $T\left(
x\right) =\left( \pi \left( x\right) ,f\left( x\right) \right) $, $x\in M^{3}
$, such that $g\circ \pi =\xi $ and $f=\Psi _{g}\circ T$. From the second
part of Proposition \ref{p5} it follows that the normal curvature of $g$ is
identically zero and the Gaussian curvature $K$ of the metric induced by $g$
satisfies $K>1$. So $g$ has not totally geodesic points. Consequently, by
virtue of Propositions \ref{p2} and\ \ref{p3}, we deduce that $g$ is the
associate of a minimal immersion $h:M^{2}\rightarrow Q^{3}$, where $Q^{3}$
is an equidistant hypersurface of $\mathbb{H}^{4}$.

\textbf{Case 2.} Assume now that the nullity distribution of $f$ doesn't
allow a global unit section. We can pick out two unit vectors $e_{2}\left(
x\right) $, $-e_{2}\left( x\right) \in \Delta \left( x\right) $, for each $%
x\in M^{3}$. Then we can construct a $2$-fold covering space $\widetilde{M}%
^{3}$ of $M^{3}$ with covering map $\Pi :\widetilde{M}^{3}\rightarrow M^{3}$
by choosing the two points in $\Pi ^{-1}\left( x\right) $ to correspond to
these vectors. One can easily check that $\widetilde{M}^{3}$ is a connected
and oriented manifold. Now we equip $\widetilde{M}^{3}$ with the covering
metric and consider the isometric immersion $\widetilde{f}:=f\circ \Pi :%
\widetilde{M}^{3}\rightarrow \mathbb{H}^{4}$ with unit normal $\widetilde{%
\xi }:=\xi \circ \Pi $. Obviously $d\Pi $ preserves the principal directions
and the principal curvatures of $\widetilde{f}$ are $\widetilde{k}_{1}=-%
\widetilde{k}_{3}=\widetilde{\lambda }:=\lambda \circ \Pi $, $\widetilde{k}%
_{2}=0$. We can readily verify that there exists a global unit vector field $%
\widetilde{e}_{2}$ which spans the nullity distribution $\widetilde{\Delta }$
of $\widetilde{f}$. It is clear that $\widetilde{f}$ satisfies all the
assumptions of Case 1. Moreover there exists a deck transformation $a:%
\widetilde{M}^{3}\rightarrow \widetilde{M}^{3}$ which is in fact an
involution and $\Pi ^{-1}\left( x\right) =\left\{ \widetilde{x},a\left( 
\widetilde{x}\right) \right\} $, for each $x\in M^{3}$.

The deck transformation $a$ induces an involution $\widetilde{a}$ on the
quotient space $\widetilde{M}^{2}$ of leaves of $\widetilde{\Delta }$ in a
natural way. Since for each $\widetilde{x}\in \widetilde{M}^{3}$ there is no
integral curve of $\widetilde{e}_{2}$ joining $\widetilde{x}$ with $a\left( 
\widetilde{x}\right) $, the involution $\widetilde{a}$ is fixed point free.
We denote by $\widetilde{\pi }:\widetilde{M}^{3}\rightarrow \widetilde{M}%
^{2} $ the quotient projection. The quotient space $\widetilde{M}^{2}/%
\widetilde{a}$ can be equipped with the structure of a 2-dimensional
manifold which makes the projection $\widetilde{\pi }_{\widetilde{a}}:%
\widetilde{M}^{2}\rightarrow \widetilde{M}^{2}/\widetilde{a}$ a covering
map. The map $q:\widetilde{M}^{2}/\widetilde{a}\rightarrow M^{2}$ given by 
\begin{equation*}
q\circ \widetilde{\pi }_{\widetilde{a}}\circ \widetilde{\pi }=\pi \circ \Pi ,
\end{equation*}
is well defined and bijection, where $M^{2}$ is the\ quotient space of
leaves of $\Delta $ and $\pi :M^{3}\rightarrow M^{2}$ is the quotient
projection. Using the map $q$ we can equip $M^{2}$ with the structure of a
2-dimensional differentiable manifold which makes $q$ a diffeomorphism and $%
\pi $ a submersion. Thus we may identify $\widetilde{M}^{2}/\widetilde{a}$
with $M^{2}$.

From Case 1, we know that the map $\widetilde{g}:\widetilde{M}%
^{2}\rightarrow \mathbb{S}_{1}^{4}$, which is induced by $\widetilde{\xi }$,
is a spacelike stationary immersion without totally geodesic points.
Furthermore, there exists a minimal immersion $\widetilde{h}:\widetilde{M}%
^{2}\rightarrow Q^{3}$, where $Q^{3}$ is an equidistant hypersurface of $%
\mathbb{H}^{4}$, such that $\widetilde{g}$ coincides with the associate of $%
\widetilde{h}$ and $\widetilde{f}=\Psi _{\widetilde{g}}\circ \widetilde{T}$,
where $\widetilde{T}$ is the isometry given by $\widetilde{T}\left( 
\widetilde{x}\right) =(\widetilde{\pi }\left( \widetilde{x}\right) ,%
\widetilde{f}\left( \widetilde{x}\right) )$, $\widetilde{x}\in \widetilde{M}%
^{3}$. Bearing in mind Propositions \ref{p2}, \ref{p3} and taking into
account $\left( 3.4\right) $ we easily see that 
\begin{equation*}
\widetilde{h}\circ \widetilde{\pi }=\frac{1}{\sqrt{1-\widetilde{v}^{2}}}%
\left( \widetilde{v}d\widetilde{f}\left( \widetilde{e}_{2}\right) +%
\widetilde{f}\right) ,
\end{equation*}
where $\widetilde{v}=\widetilde{e}_{2}(\log \widetilde{\lambda })$. Since $a$
is a deck transformation, the maps $g:M^{2}\rightarrow \mathbb{S}_{1}^{4}$
and $h:M^{2}\rightarrow Q^{3}$ given by 
\begin{equation*}
g\circ \widetilde{\pi }_{\widetilde{a}}=\widetilde{g}\text{ \ and \ }h\circ 
\widetilde{\pi }_{\widetilde{a}}=\widetilde{h},
\end{equation*}
are well defined. Then $g$ is the associate of $h$, since $\widetilde{g}$ is
the associate of $\widetilde{h}$. Moreover $f=\Psi _{\widehat{h}}\circ T$,
where $T:M^{3}\rightarrow \mathcal{N}^{1}(\widehat{h})$ is the local
isometry given by $T\left( x\right) =\left( \pi \left( x\right) ,f\left(
x\right) \right) $, $x\in M^{3}$. This completes the proof.
\end{proof}

\begin{remark}
It should be interesting to know whether a similar classification result can
be obtained without the assumption that the scalar curvature is bounded from
below.
\end{remark}

\end{document}